\documentclass[12pt, reqno]{amsart}
\usepackage{amssymb}
\usepackage{amsmath}
\usepackage{bbm}
\usepackage{amsthm}
\usepackage{mathtools}
\usepackage{graphicx}
\usepackage{amscd}
\usepackage{epic}
\usepackage[mathscr]{euscript}
\usepackage{comment}

\usepackage{geometry}
\usepackage{bm}

\usepackage[dvipsnames]{xcolor}
\usepackage[colorlinks,citecolor=OliveGreen,linkcolor=Mahogany,urlcolor=Plum]{hyperref}
\usepackage[alphabetic]{amsrefs}

\usepackage{todonotes}

\theoremstyle{plain}
\newtheorem{theorem}{Theorem}[section]
\newtheorem{lemma}[theorem]{Lemma}
\newtheorem{proposition}[theorem]{Proposition}

\newtheorem{thm}[theorem]{Theorem}

\theoremstyle{definition}

\newtheorem{definition}[theorem]{Definition}

\theoremstyle{remark}
\newtheorem{remark}[theorem]{Remark}

\newcommand{\cB}{\mathcal{B}}

\newcommand{\cF}{\mathcal{F}}

\newcommand{\cO}{\mathcal{O}}

\newcommand{\cS}{\mathcal{S}}

\newcommand{\cW}{\mathcal{W}}

\newcommand{\bC}{\mathbb{C}}

\newcommand{\bP}{\mathbb{P}}
\newcommand{\bQ}{\mathbb{Q}}
\newcommand{\bR}{\mathbb{R}}

\newcommand{\bW}{\mathbb{W}}
\newcommand{\bZ}{\mathbb{Z}}

\newcommand{\fa}{\mathfrak{a}}
\newcommand{\fb}{\mathfrak{b}}

\newcommand{\fm}{\mathfrak{m}}

\newcommand{\mult}{\mathrm{mult}}

\newcommand{\vol}{\mathrm{vol}}
\newcommand{\ord}{\mathrm{ord}}
\newcommand{\lct}{\mathrm{lct}}

\newcommand{\Val}{\mathrm{Val}}

\newcommand{\Spec}{\mathrm{Spec~}}

\newcommand{\mld}{\mathrm{mld}}
\newcommand{\gr}{\mathrm{gr}}

\newcommand{\im}{\mathrm{im}}

\newcommand{\e}{\mathrm{e}}

\def\DH{\mathrm{DH}}

\def\dif{{\rm d}}

\newcommand{\bu}{\bullet}

\newcommand{\nvol}{\widehat{\mathrm{vol}}}
\newcommand{\NVOL}{\widehat{\mathrm{Vol}}}

\newcommand{\hvol}{\widehat{\mathrm{vol}}}

\newcommand{\BH}{\mathbf{H}^\mathrm{NA}}

\numberwithin{equation}{section}

\begin{document}

\title{Stable degeneration and birational geometry}

\author{Lu Qi}\footnote{The author is partially supported by a Shanghai Sailing program 24YF2709800.}
\address{School of Mathematical Sciences, East China Normal University, Shanghai 200241, China}
\email{lqi@math.ecnu.edu.cn}

\begin{abstract}
     This expository article is based on the author’s talk at the Kinosaki Algebraic Geometry Symposium 2025. We discuss some recent progress surrounding stable degeneration in algebraic K-stability theory. 
\end{abstract}

\maketitle

\section{Introduction}

\emph{Throughout, we work over the field $\bC$ of complex numbers.}

\medskip

The local (algebraic) theory of K-stability for Kawamata log terminal (klt) singularities was initiated in \cite{Li-nv}.  C. Li introduced the \emph{normalized volume} functional 
\begin{align*}
    \nvol_{(X,\Delta),x}:\Val_{X,x}&\to \bR_{>0}\cup\{+\infty\}\\
        v&\mapsto\left\{\begin{aligned}
            &A_{X,\Delta}(v)^n\cdot \vol(v), &\text{if}~ A_{X,\Delta}(v)<+\infty,\\
            &+\infty, &\text{otherwise.}
        \end{aligned}\right.
\end{align*}
The \emph{local volume} of the singularity is then defined to be
\begin{equation}\label{eqn:local vol}
    \nvol(x,X,\Delta)=\inf_{v\in\Val_{X,x}} \nvol_{(X,\Delta),x}(v)
\end{equation}

One central topic is minimizing the local volume functional. The Stable Degeneration Conjecture \cites{Li-nv,LX-KC} predicts that the minimization process gives rise to a two-step degeneration, which is now understood by intensive work including \cites{Blu-existence,LX18,LWX18,Xu-qm,XZ-uniqueness,XZ-sdsing} (see also \cites{LX-KC,BLQ-convexity}).

For simplicity, we work with $\bQ$-pairs, i.e., pairs whose coefficients belong to $\bQ\cap (0,1)$. 
See \cites{HLQ-vol-ACC,Zhu-mld^K-2} for the generalization to $\bR$-pairs.

\begin{theorem}[Stable degeneration for singularities]\label{thm:stable degeneration for klt singularities}
    Let $x\in (X,\Delta)$ be a klt singularity, where $X=\Spec R$. Then the following statements hold.
    \begin{enumerate}
        \item (Existence) There exists a valuation $v_*\in\Val_{X,x}$ such that 
        \[
            \nvol(x,X,\Delta)=\nvol_{(X,\Delta),x}(v_*).
        \]
    
        \item (Uniqueness) The minimizer $v_*$ is unique up to rescaling; that is, for any valuation $v'\in\Val_{X,x}$ satisfying $\nvol(x,X,\Delta)=\nvol_{(X,\Delta),x}(v')$, there exists $\lambda\in\bR_{>0}$ such that $v'=\lambda v_*$. 

        \item (Quasi-monomiality) The minimizer $v_*$ is quasi-monomial.
    
        \item (Finite generation) The associated graded ring $\gr_{v_*} R$ is finitely generated. 
    
        \item (Two-step degeneration) The minimizer $v_*$ induces a special degeneration 
        \[
            x\in (X,\Delta)\leadsto x_0\in (X_0,\Delta_0,\xi_0)
        \]
        to a K-semistable log Fano cone singularity, where $X_0\coloneqq \Spec \gr_{v_*} R$ and $\xi_0$ is the Reeb vector induced by $v_*$.      
        Furthermore, there is a unique special degeneration 
        \[
            x_0\in (X_0,\Delta_0)\leadsto x_p\in (X_p,\Delta_p,\xi_0)
        \]
        to a K-polystable log Fano cone singularity.
    \end{enumerate}
\end{theorem}

We refer readers to \cite{LLX18} for a guide to the area and to \cite{Zhu-localsurvey} for an excellent, more comprehensive survey of recent developments. 


\medskip

Next, let us outline the proof of stable degeneration for klt singularities.

\textbf{Existence}. 
The existence of a $\nvol$-minimizer is first proven in \cite{Blu-existence}.
To find a critical point of a functional $F$, possibly enlarging the domain of definition, one first finds a sequence such that the values of $F$ approach the infimum. 
One then shows that the approximation sequence admits a limit (compactness), which lies in the original space (regularity) such that the infimum is actually achieved (continuity).

For the approximation step, a key observation of \cite{Liu-vol-sing-KE} is that the infimum in \eqref{eqn:local vol} is equal to the infimum of the \emph{normalized multiplicity} of all \emph{graded sequences of $\fm$-primary ideals} or filtrations $\fa_\bullet$, that is,
\begin{equation}\label{eqn:normalized mult}
    \nvol(x,X,\Delta)=\inf_{\fa_\bullet} \lct(X,\Delta;\fa_\bullet)^n\cdot \mult(\fa_\bullet).
\end{equation}

The existence of a limit point of the approximating sequence is achieved by the \emph{generic limit} construction; see, for example, \cites{dFM-generic-limit,Kol-which-power}. 
The idea is to parameterize all graded sequences $\fa_\bullet$ whose normalized multiplicity is bounded from above. The compactness follows from Chevalley's theorem for constructible sets, and to achieve continuity, one needs some uniform control for the convergence rate of the volume.

\textbf{Uniqueness}. 
The uniqueness of the $\nvol$-minimizer is first proven in \cite{XZ-uniqueness}, where K-stability for a valuation on a klt singularity is introduced. The proof relies implicitly on an idea of convexity.  

Later, the idea is explicitly realized in another proof given in \cite{BLQ-convexity}, where the convexity of multiplicities along \emph{geodesics} is shown. In general, the linear combination of two valuations is not a valuation. This is where \eqref{eqn:normalized mult} shows its power again. For two filtrations $\fa_{\bullet,0}$ and $\fa_{\bullet,1}$, the \emph{geodesic} $\fa_{\bullet,t}$ is defined to be the filtration
\[
    \fa_{\lambda,t}\coloneqq \sum_{(1-t)\mu+t\nu=\lambda} \fa_{\mu,0}\cap \fa_{\nu,1}.
\]

As the term suggests, the geodesic plays the role of linear combinations in the space of (saturated) filtrations. The upshot of \cite{BLQ-convexity} is that the function
\[
    t\mapsto \mult(\fa_{\bullet,t})
\]
is strictly convex along the geodesic unless the endpoints differ by a rescaling. 
Combined with \cite[Theorem 3.11]{XZ-uniqueness}, which implies the convexity of the function
\[
    t\mapsto \lct(X,\Delta;\fa_{\bullet,t}),
\]
one immediately obtains uniqueness by \eqref{eqn:normalized mult}.

\textbf{Quasi-monomiality}. 
The statement that any $\nvol$-minimizer is quasi-monomial was proven in \cite{Xu-qm}, where the argument can be divided into two steps, as in the proof for existence. 

The approximation step is the following result of \cite{LX-KC},
\[
    \nvol(x,X,\Delta)=\inf_S \nvol_{(X,\Delta),x}(\ord_S),
\]
where the infimum is now taken over all \emph{Koll\'ar components}, which means a prime divisor $S$ on a proper birational model $\mu: Y\to (X,\Delta)$ such that $-S$ is a $\bQ$-Cartier $\mu$-ample divisor and that $(Y,\mu^{-1}_*\Delta+S)$ is plt. 

Compactness follows from a combination of the \emph{boundedness of complements} and a generic limit argument.
Indeed, it is not hard to see that any Koll\'ar component $S$ is a log canonical place of some $\bQ$-complement and hence of some $N$-complement for some $N$ depending only on $X$, by \cite[Theorem 1.8]{Birkar-bab-1}. 
However, unlike the projective setting, one is still not done since the above statement does not yield a finite-dimensional, thus compact, parameter space. 
The solution is to show, using a properness estimate, that all information can be extracted up to the quotient by a fixed power of the maximal ideal. 

Note that the argument yields another proof of the existence of the minimizer.


\textbf{Finite generation}. 
The finite generation problem turns out to be the most subtle part of the stable degeneration conjecture.
In the \emph{quasi-regular} case, namely when the minimizer $v_*$ is divisorial, finite generation is known by \cites{LX-KC,Blu-existence}, which follows essentially from \cite{BCHM10}. 
In dimension $2$, it is known by \cites{Cut-fg} (see also \cite{LLX18}), essentially because a non-divisorial minimizer in this case is toric.
The finite generation for a general higher rank valuation is settled in \cite{XZ-sdsing}, which is a local analog of the theory called \emph{higher rank finite generation} established in \cite{LXZ-HRFG}.
The strategy is to show that the minimizer of the functional is a \emph{special} valuation (which is quasi-monomial a priori), whose associated graded ring is then shown to be finitely generated.

Recently, \cite{Che-HRFG} has provided an elegant and more algebraic approach to the finite generation problem, which works in a slightly more general setting.

\section{Applications}

Next, we move on to a brief discussion of some applications of stable degeneration for singularities. 

\subsection{Ingredients and corollaries}

First, we discuss two byproducts of the proofs in the last section, as well as an immediate corollary of the stable degeneration theorem.

\medskip

\subsubsection{Complex geometry}
During the course of the proof for quasi-monomiality in \cite{Xu-qm}, Xu actually proves the following stronger statement, which verifies the weak version of \cite[Conjecture B]{JM12} and completes the algebraic approach proposed in \emph{op. cit.} toward the openness conjecture of \cite{DK-openness}.

\begin{theorem}\label{thm:JM conjecture}
    Let $(X,\Delta)$ be a klt log pair. For any filtration $\fa_\bullet$ on $X$ such that $\lct(X,\Delta;\fa_\bullet)<\infty$, there exists a quasi-monomial valuation $v_*$ such that 
    \[
        \lct(X,\Delta;\fa_\bullet)=\frac{A_{X,\Delta}(v_*)}{v_*(\fa_\bullet)}
    \]
\end{theorem}

The quasi-monomiality of $\nvol$-minimizers follows from the above theorem, since one can show that any $\nvol$-minimizer $v_*$ is the unique valuation (up to scaling) that computes $\lct(X,\Delta;\fa_\bullet(v_*))$.

\subsubsection{Commutative algebra}
In order to formulate the condition for the linearity of multiplicities, \cite{BLQ-convexity} introduces the notion of \emph{saturation} $\widetilde{\fa_\bullet}$ for a filtration $\fa_\bullet$ on a Noetherian local domain.  
As a side product, one can show the following generalization of a classical result of Rees \cite{Ree-mult-ideal} that the Hilbert-Samuel multiplicities of two $\fm$-primary ideals $\fa\subset\fb$ are equal if and only if they have the same integral closure. 

\begin{theorem}\label{thm:Rees}
    Let $(R,\fm)$ be a Noetherian local domain that is analytically irreducible. Let $\fa_\bullet\subset\fb_\bullet$ be two $\fm$-filtrations. Then $\e(\fa_\bullet)=\e(\fb_\bullet)$ if and only if $\widetilde{\fa_\bullet}=\widetilde{\fb_\bullet}$. 
\end{theorem}

Saturation also provides a notion to formulate when the Minkowski equality for filtrations holds; see, for example, \cite{Cut-Minkowski-eq}. 

\subsubsection{Moduli theory for Fano varieties}
As we have seen, the local theory relies on techniques inspired by the proofs of the global results. In return, the local results, especially the following corollary of uniqueness proven in \cite{XZ-uniqueness}, have applications in global geometry. 

\begin{theorem}\label{thm:bdd of index}
    Let $x\in (X,\Delta)$ be a klt singularity. 
    There exists $r\in\bZ_{>0}$ such that $r D$ is Cartier for any $\bQ$-Cartier integral divisor $D$.
\end{theorem}

The index control can be powerful on many occasions and plays an important role in the so-called \emph{moduli continuity method} in the theory of K-moduli for Fano varieties; see, for example, \cites{LX19,SS-GH-compact,ADL21,Liu-cubic-4,LZ-genus-4}.

We remark that Theorem \ref{thm:bdd of index} is often applied in conjunction with the local-global volume comparison explored in \cites{Fuj18,Liu-vol-sing-KE}. 
Shortly, we will see how they can be applied to study the termination of the Minimal Model Program. 

\subsection{Applications to birational geometry}

Since \cite{LX14}, birational geometry, especially the Minimal Model Program (MMP), has been a crucial tool in the theory of K-stability. More recently, it was found in \cite{HQZ25} that K-stability can be applied to study birational geometry, the main results of which we outline below.


\subsubsection{Local-global principle for volumes}

In view of Theorem \ref{thm:bdd of index}, one main reason we can apply K-stability to study the local birational invariants is the following local-global principle for volumes, which is inspired by \cite[Theorem D]{BJ20} in the projective setting and \cite[Lemma 2.13]{XZ-local-bdd} in the local setting. 
Recall that the \emph{alpha invariant} \cite{Tia87} is defined by
\[
    \alpha(V_\bullet)\coloneqq \inf\{\lct(X,\Delta;\frac{1}{m}D)\mid D\in V_m\}.
\]

\begin{lemma}\label{lem:local-global volume}
    Let $(X,\Delta)$ be a klt log pair and let $L$ be a $\bQ$-Cartier $\bQ$-divisor on $X$. Then for any closed point $x\in X$ and any eventually birational linear series $V_\bullet$ of $L$, we have
    \[
        \nvol(x,X,\Delta)\ge \alpha(X,\Delta;V_\bullet)^n\cdot \vol(V_\bullet).
    \]
\end{lemma}

Here, a graded linear series $V_\bullet\subset \oplus_m H^0(X,mL)$ is \emph{eventually birational} if and only if for sufficiently divisible $m$, the rational map $\phi_{V_m}:X\to \bP(V_m^\vee)$ is birational onto its image. 
For example, the complete linear series of a big line bundle on a projective variety is eventually birational. 
This assumption, together with the condition that $\alpha(V_\bullet)>0$, is enough to guarantee that $V_\bullet$ has positive \emph{volume}
\[
    \vol(V_\bullet)\coloneqq\lim_{m\to\infty}\frac{\dim(V_m)}{m^n/n!}.
\]

Now we can define the log canonical volume of a graded linear series. For simplicity, we leave out some technical details and only consider log pairs here.

\begin{definition}\label{defn:log canonical volume}
    Let $(X,\Delta)$ be a klt log pair and let $L$ be a $\bQ$-divisor on $X$. The \emph{log canonical volume} of $L$ is defined to be
    \[
        \NVOL_{X,\Delta}(L)\coloneqq\sup_{\alpha(X,\Delta;V_\bullet)\ge 1}\{\vol(V_\bullet)\}.
    \]
\end{definition}

A birational contraction $\phi:(X,\Delta)\dashrightarrow (X',\Delta')$ is called an \emph{MMP-type contraction} if there exist proper birational morphisms $p:W\to X$ and $q:W\to X'$ such that $g=\phi\circ p$ and that
\[
    p^*(K_X+\Delta)-q^*(K_{X'}+\Delta')\ge 0.
\]
Intuitively, an MMP-type contraction does not worsen the singularities, and hence $\NVOL$ is non-decreasing along MMP-type contractions almost by definition. By contrast, the local volume $\nvol$ does not seem to satisfy any monotonicity.

\begin{lemma}\label{lem:positive lc volume}
    Let $(X,\Delta)$ be a klt log pair and let $L\ge 0$ be a big $\bQ$-divisor such that $(X,\Delta+L)$ is klt. Then $\NVOL_{X,\Delta}(L)>0$. 
\end{lemma}

Recall that general type MMPs can be reduced to the big boundary case by a standard argument (for example, \cite[Lemma 3.7.5]{BCHM10}). Combining Lemma \ref{lem:local-global volume} and Lemma \ref{lem:positive lc volume}, we know that the local volume stays away from $0$ in such MMPs. 

\begin{thm}\label{thm:vol bdd}
    Let $(X,\Delta)$ be a projective $\bQ$-factorial klt pair such that $\Delta$ is big. Then there exists some $\varepsilon>0$ such that for any sequence of a $(K_X+\Delta)$-MMP $(X,\Delta)\dashrightarrow (X',\Delta')$ and any closed point $x'\in X'$, we have $\hvol(x',X')\ge \varepsilon$.
\end{thm}

\subsubsection{Boundedness in general type MMP}

In view of Theorem \ref{thm:bdd of index}, the lower bound for local volumes of Theorem \ref{thm:vol bdd} immediately gives a uniform bound for the Cartier index in any general type MMP, which implies the discreteness of minimal log discrepancies.
Moreover, applying \cite[Theorem 1.3]{XZ-local-bdd}, we also know that the minimal log discrepancy is uniformly bounded from above. Hence, we obtain the following result, fitting into Shokurov's approach \cite{Sho-mld-conj} toward the termination of flips. 

\begin{thm}\label{thm:index+mld}
    Let $(X,\Delta)$ be a projective $\bQ$-factorial klt pair such that $\Delta$ is big. Then there exist $r\in\bZ_{>0}$ and a finite set $S\subseteq \bR_{>0}$, depending only on the pair $(X,\Delta)$, such that for any sequence of steps of a $(K_X+\Delta)$-MMP $(X,\Delta)\dashrightarrow (X',\Delta')$,
    \begin{enumerate}
        \item $rD'$ is Cartier for any Weil divisor $D'$ on $X'$, and
        \item for any point $x'\in X'$, we have $\mld(x',X',\Delta')\in S$.
    \end{enumerate}
\end{thm}

Another strategy toward termination is the finiteness of models, which has been employed in \cite{BCHM10}. In this direction, we can prove the following result.

\begin{proposition}\label{prop:boundedness implies termination}
    Let $(X_0,\Delta_0) \coloneqq (X,\Delta)$ be a log canonical pair and let 
    \[
        (X_0,\Delta_0)\dashrightarrow (X_1,\Delta_1) \dashrightarrow \cdots \dashrightarrow (X_i,\Delta_i) \dashrightarrow  \cdots
    \]
    be a sequence of a $(K_X+\Delta)$-MMP. Assume that the set of pairs $\{(X_i,\Delta_i)\,|\,i\in\bZ_{\ge 0}\}$ is log bounded. Then this MMP terminates after finitely many steps. 
\end{proposition}

So far, we cannot bound all models appearing in an MMP, but the tools from K-stability allow us to bound the fibers of all the contractions. 

\begin{thm}\label{thm:fiber}
    Let $(X,\Delta)$ be a projective $\bQ$-factorial klt pair such that $\Delta$ is big. Then there exists a projective family $\cW\to \cB$ over a finite type base $\cB$, such that in any sequence of $(K_X+\Delta)$-MMP, every fiber of the extremal contractions or the flips is isomorphic to $\cW_b$ for some $b\in \cB$.
\end{thm}

There are two main ingredients for the above theorem.
The first one is a boundedness type result in the spirit of \cite[Theorem 1.1]{XZ-local-bdd} (see also \cite{Jia17} for the global version), whose proof indeed relies heavily on \emph{loc. cit.} and the relative cone construction. Roughly speaking, log Fano fibrations with $\NVOL$ bounded away from $0$ form a specially bounded family.
The other ingredient is that special boundedness implies the boundedness of fibers, which is inspired by the proof of \cite[Theorem 1.3]{XZ-local-bdd}. Here, the strategy is somewhat brutal. Via a syzygy argument, the boundedness of the central fiber of the special degeneration controls the embedded dimension as well as the degree of the defining equations for any fiber of the fibration.

\section{Generalizations and discussions}

\subsection{Stable degeneration for Fano fibration germs}

We next turn to a generalization of the stable degeneration theory to the setting of Fano fibration germs, for which a similar $2$-step degeneration framework was conjectured in \cite{SZ24}.
A \emph{(log) Fano fibration} $f:(X,\Delta)\to Z$ is a surjective projective morphism between normal varieties such that 
\begin{enumerate}   
    \item[(0)] $f$ is a fibration, i.e., $f_*\cO_X=\cO_Z$, 
    
    \item $(X,\Delta)$ has klt singularities, and 

    \item $-(K_X+\Delta)$ is ample over $\bZ$. 
\end{enumerate}
A \emph{Fano fibration germ} consists of a Fano fibration $f:(X,\Delta)\to Z$, where $Z$ is affine, together with a fixed closed point $o\in Z$.

As one easily sees, a projective klt pair $(X,\Delta)$ is automatically a Fano fibration germ over $\Spec \bC$, and a klt singularity $x\in (X,\Delta)$ (where we assume that $X$ is affine by convention) is a Fano fibration germ with respect to the identity morphism. Therefore, Fano fibration germs can be viewed as a natural interpolation between klt Fano pairs and klt singularities.
Moreover, for any $(K_X+\Delta)$-negative contraction $(X,\Delta)\to Z$ and any closed point $z\in Z$, we get a Fano fibration germ $(X,\Delta)\to Z\ni z$ after replacing $Z$ with an affine neighborhood of $z$. 

In this subsection, we will work with a fixed Fano fibration germ $(X,\Delta)\to Z\ni o$ with $Z=\Spec R_0$ and $L\coloneqq -(K_X+\Delta)$. Fix $r\in\bZ_{>0}$ such that $rL$ is Cartier and very ample over $Z$. Let $R\coloneqq \oplus_{m\in r\bZ_{>0}} R_m$, where $R_m\coloneqq H^0(X,mL)$ is a finite $R_0$-module. Write $n=\dim X$ and $\fm\subset R_0$ for the maximal ideal of $o\in Z$.

In parallel with the known cases, the stable degeneration theory for Fano fibration germs is centered around the minimization of a canonical functional on the space of filtrations, the non-archimedean $\mathbf H$-invariant, which we shall proceed to define now. 

A \emph{linearly bounded $\fm$-filtration} $\cF$ on $R=\oplus_m R_m$ is a filtration $\{\cF^\lambda R_m\}_{\lambda\in\bR}$ for each $m\in r\bZ_{>0}$ satisfying the following conditions:
\begin{enumerate}
    \item (decreasing) $\cF^\lambda R_m \subset \cF^{\lambda'} R_m $ if $\lambda > \lambda'$,
        
    \item (left continuous) $\cF^\lambda R_m=\cF^{\lambda-\epsilon}R_m$ for any $0<\epsilon \ll 1$, 
        
    \item (multiplicative) $\cF^\lambda R_m \cdot \cF^{\lambda'}R_{m'} \subset \cF^{\lambda+\lambda'}R_{m+m'}$,

    

    \item (supported over $\fm$) $\cF^\lambda R_0\subset R_0$ is an $\fm$-primary ideal for any $\lambda>0$.

    \item (left linearly bounded) there exist $c\in\bR_{>0}$ and $e_-\in\bR$ such that 
    \[
        \fm^{\lceil\frac{\lambda - me_-}{c} \rceil} R_m \subset \cF^\lambda R_m 
    \]
    for any $\lambda$, and
    
    \item (right linearly bounded) there exist $C\in\bR_{>0}$ and $e_+\in\bR$ such that 
    \[
    \cF^\lambda R_m \subseteq \fm^{\lceil\frac{\lambda - me_+}{C} \rceil} R_m
    \]
    for any $\lambda$. 
\end{enumerate}

A typical example is the filtration $\cF_v$ induced by a valuation $v$ on $X$, defined by
\[
    \cF^\lambda_vR_m\coloneqq \{s\in R_m\mid v(s)\ge \lambda\}.
\]
Note that $\cF_v$ satisfies (1)-(5) above if the center of $v$ is contained in $f^{-1}(o)$, and it satisfies (6) if $A_{X,\Delta}(v)<+\infty$, by the Izumi inequality \cite[Theorem 1.2]{Li-nv}. In what follows, a filtration on $R$ means a linearly bounded $\fm$-filtration. 

\medskip

As in the case of normalized volumes, the $\BH$-invariant consists of a term encoding the singularities and a term containing the information regarding positivity. 

\begin{definition}\label{defn:log canonical slope}
Let $\cF$ be a filtration on $R$. For $\lambda\in\bR_{>0}$, the graded sequence of {\it base ideals} $I^{(\lambda)}_\bullet=I^{(\lambda)}_\bullet(\cF) = \{I_{m,m\lambda}\}_{m\in \bZ_{>0}}$ of $\cF$ is defined by 
\begin{equation*}
\label{eqn:lc slope}
I_{m,m\lambda}=
I_{m,m\lambda}(\cF)
\coloneqq
\im(\cF^{m\lambda}R_m\otimes \cO_X(-mL)\to \cO_X)
\end{equation*}
for any $m\in r\bZ_{>0}$. 
\end{definition}

Note that given a filtration $\cF$, $\gr^\lambda_\cF R_m\coloneqq \cF^\lambda R_m/\cF^{>\lambda} R_m$ is a finite-dimensional vector space over $\bC$. So we can define a sequence of Duistermaat-Heckman measures.

\begin{definition}\label{defn:discrete DH measure}
    Let $\cF$ be a filtration on $R$. For $m\in r\bZ_{>0}$, define
    \[
        \DH_{\cF,m}\coloneqq \frac{1}{m^n/n!}\sum_{\lambda\in\bR} \dim(\gr^{m\lambda}_\cF R_m)\cdot \delta_\lambda=\frac{1}{m^n/n!}\frac{\dif}{\dif \lambda}\dim (R_m/\cF^{m\lambda} R_m).
    \]
\end{definition}

The following result essentially follows from the argument of \cite[Appendix A]{SZ24}. 

\begin{proposition}\label{prop:DH measure}
    The measures $\DH_{\cF,m}$ converge weakly to the measure
    \[
        \DH_\cF\coloneqq \frac{\dif}{\dif \lambda} \vol(R/\cF^{(\lambda)}R),
    \]
    where 
    \[
        \vol(R/\cF^{(\lambda)}R)\coloneqq \limsup_{m\to\infty} \frac{\dim (R_m/\cF^{m\lambda} R_m)}{m^n/n!}.
    \]
\end{proposition}

Now we are ready to define the $\BH$-invariant.

\begin{definition}\label{defn:invariants}
    Let $\cF$ be a filtration on $R$. 
    \begin{enumerate}
        \item The \emph{log canonical slope} of $\cF$ is defined by
            \[
                \mu(\cF)=\mu_{X,\Delta}(\cF)\coloneqq \sup\{\lambda\mid \lct(X,\Delta;I^{(\lambda)}_\bu)\ge 1\}. 
            \]

        \item The \emph{$\widetilde S$-invariant} of $\cF$ is defined by
        \[
            \widetilde S(\cF)\coloneqq -\log \int_\bR e^{-\lambda}\DH_\cF(\dif\lambda).
        \]

        \item Define
        $\BH(\cF)\coloneqq \mu(\cF)-\widetilde S(\cF)$, and 
        \[
            h((X,\Delta)\to Z\ni o)\coloneqq \inf_\cF \BH(\cF),
        \]
        where the infimum is taken over all linearly bounded $\fm$-filtrations $\cF$. 
    \end{enumerate}
\end{definition}

The following stable degeneration theorem for Fano fibration germs, parallel to the local case, will be proven in \cite{HMQWZ-sdFFG}, which verifies \cite[Conjecture 6.4]{SZ24}.

\begin{theorem}[Stable degeneration for Fano fibration germs]\label{thm:stable degeneration for Fano fibration germs}
Let $(X,\Delta)\to Z\ni o$ be a log Fano fibration germ. Then the following statements hold.
\begin{enumerate}
    \item There exists a valuation $v_*\in\Val_{X,o}$ such that 
    \[
        \BH(v_*) = h((X,\Delta)\to Z\ni o).
    \]
    
    \item The minimizer $v_*$ is unique. 

    \item The minimizer $v_*$ is quasi-monomial.
    
    \item The associated graded ring $\gr_{v_*} R$ is finitely generated. 
    
    \item The minimizer $v_*$ induces a special degeneration 
        \[
            (X,\Delta)\to Z\ni o \quad\leadsto\quad(X_{0},\Delta_{0},\xi_0)\to Z_{0}  \ni o
        \]
        to a K-semistable polarized Fano fibration germ.      
        Furthermore, there is a unique special degeneration 
        \[
            (X_{0},\Delta_{0},\xi_0)\to Z_{0}  \ni o\quad \leadsto\quad (X_{p},\Delta_{p},\xi_0)\to Z_{p}\ni o.
        \]
\end{enumerate}
\end{theorem}

\begin{remark}
    \begin{enumerate}
        \item The conjecture in \cite{SZ24} was stated in terms of the \emph{weighted volume} functional $\bW:\Val_{X,o}\to \bR_{>0}$, which is defined by
        \[
            \bW(v)\coloneqq e^{A_{X,\Delta}(v)}\int_\bR e^{-\lambda} \DH_v(\dif \lambda).
        \]
        It can be shown that if $v$ is \emph{weakly special}, i.e., an lc place of a complement, then $A_{X,\Delta}(v)=\mu_{X,\Delta}(\cF_v)$. Moreover, there exists a sequence of weakly special valuations $v_i$ such that $\lim_i\BH(\cF_{v_i})=h((X,\Delta)\to Z\ni o)$. 
        Therefore, we have a Liu-type regularization
        \[
            h((X,\Delta)\to Z\ni o)=\inf_\cF \BH(\cF)=\inf_v \log\bW(v),
        \]
        which reduces the original conjecture to the above theorem.

        \item In the local case where $Z=X=\Spec R_0$ and $x\in (X,\Delta)$ is a klt singularity, a filtration $\cF$ on $R$ is the same as a filtration $\fa_\bullet$ on $R_0$. It is not hard to show that $\mu(\cF)=\lct(\fa_\bullet)$ and that $\widetilde S(\cF)=\log \mult(\fa_\bullet)$ (the latter is essentially observed in \cite[Appendix A]{BLQ-convexity}).
        Hence 
        \[
            e^{\BH(v)}=e^{A_{X,\Delta}(v)}\vol(v)
        \]
        for weakly special valuations $v$, and it is a calculus exercise to check that the minimizer for $\BH$ is exactly the minimizer of $\nvol$ satisfying $A(v)=n$. This explains why the uniqueness up to scaling in the local case becomes genuine uniqueness in our setting; see \cite[Example 4.4]{SZ24}.

        \item A formula for $\bW$ in terms of restricted volumes can be found in \cite{Oda-SZ}, together with several concrete examples and estimates.
    \end{enumerate}
\end{remark}




\subsection{Boundedness results and comparison of the invariants}

A vital component of the K-stability theory concerns \emph{boundedness} in a suitable sense. A major step in the proof of Theorem \ref{thm:fiber} is the following statement.

\begin{theorem}\label{thm:special boundedness}
    Fix $n\in \bZ_{>0}$, $\varepsilon\in\bR_{>0}$, and let $I\subseteq [0,1]\cap \bQ$ be a finite set. Then the set 
    \begin{align*}
        \cS\coloneqq\left\{
            \begin{aligned}
                &(X,\Delta)\to Z\ni z~\text{is a}\\
                &\text{Fano fibration germ}
            \end{aligned}
            \left\vert
            \begin{aligned}
            &\dim X=n,~ \Delta\in I,~\NVOL_{X,\Delta}(H)\ge \varepsilon\\
            &\text{for some $\bQ$-divisor $H\sim_{\bQ,Z} -(K_X+\Delta)$}
            \end{aligned}\right.\right\}
    \end{align*}
    is \emph{specially bounded}. 
\end{theorem}

We skip the precise definition of special boundedness here, but the statement eventually boils down to a boundedness result for Fano fibration germs with a torus action whose log canonical volumes are bounded away from $0$. 
As noted earlier, the log canonical volume is an intermediate invariant used to control the local volume, while the $\BH$-invariant (or the weighted volume) is an invariant carrying certain differential-geometric information. 

Therefore, the author is tempted to ask whether these invariants are actually comparable, and thus one can replace $\NVOL$ with $\BH$ in Theorem \ref{thm:special boundedness} as in the local case.


\subsection*{Acknowledgments} The author would like to thank Masafumi Hattori and the organizers of the Kinosaki Algebraic Geometry Symposium 2025, Makoto Enokizono, Kenta Hashizume, JuAe Song, and Kazuhiko Yamaki, for their kind invitation and hospitality. 
The author would also thanks Ziquan Zhuang for kind comments on a preliminary version of this article.

\bibliography{ref}

\end{document}